\documentclass[a4paper,12pt]{elsarticle}
\usepackage{graphicx}
\usepackage[all]{xy}
\usepackage[english]{babel}
\usepackage{amssymb,amsmath,amsthm}

\usepackage{pstricks,pst-plot}
\usepackage{verbatim}
\usepackage{hyperref}

\newtheorem{theorem}{Theorem}[section]
\newtheorem{question}{Problem}
\newtheorem*{proba}{Problem 1a}

\newtheorem{lemma}[theorem]{Lemma}
\newtheorem{proposition}[theorem]{Proposition}

\newtheorem{corollary}[theorem]{Corollary}

\newdefinition{definition}{Definition}
\newdefinition{example}{Example}
\newdefinition{remark}{Remark}[theorem]
\newcommand{\ff}{\mathbb F}

\newcommand{\fq}{\mathbb{F}_q}
\newcommand{\fqn}{\mathbb{F}_{q^n}}
\newcommand{\fqd}{\mathbb{F}_{q^d}}
\newcommand{\F}{\mathbb{F}}

\newcommand{\M}{\mathcal W}
\newcommand{\V}{\mathcal V}
\newcommand{\rc}{\mathcal R}
\newcommand{\mvsp}{\operatorname{m.v.s.p.}}
\newcommand{\mvsps}{\operatorname{m.v.s.p.}\text{'s}}

\newcommand{\gal}{\operatorname{Gal}}

\newcommand{\Conceicao}{Concei{\c c}\~ao}

\title{On the characterization of minimal value set polynomials}
\begin{document}
\begin{frontmatter}

%\email{rconceic@math.utexas.edu}
%\subjclass[2000]{Primary 54C40, 14E20; Secondary 46E25, 20C20}
\title{On the characterization of minimal value set polynomials}

\author[herivas]{Herivelto Borges}
\author[rico]{Ricardo \Conceicao}
\address[herivas]{Universidade de S\~ao Paulo, Inst. de Ci\^encias Matem\'aticas e de Computa\c c\~ao, S\~ao Carlos, SP 13560-970, Brazil.}
\address[rico]{Oxford College of Emory University. 100 Hamill Street, Oxford, Georgia 30054.}

\date{July 2010}

\begin{abstract} 
%In this paper,
 We give an explicit characterization of all minimal value set polynomials in $\F_q[x]$ whose  set of values is a subfield $\F_{q'}$ of $\F_{q}$. We show that the set of such polynomials, together with the constants of $\F_{q'}$, is an  $\F_{q'}$-vector space of dimension $2^{[\F_{q}:\F_{q'}]}$. Our approach   not only provides the exact number of such polynomials, but also yields a construction  of new examples of minimal value set polynomials for some other fixed value sets. In the latter case, we  also derive a non-trivial  lower bound for the number of such polynomials.

\end{abstract}


\end{frontmatter}

%\tableofcontents
\section{Introduction}\label{sec:intro}

%\textcolor{red}{Uma ``rough idea" de como eu acho que a introducao e as preliminares devem ser.}

%If the polynomials $f$ and $g$ have the same value set $V$, then the curve $f(x)=g(y)$ is more likely to have many points when $\mid V \mid$ is small(see Proposition \ref{lagrange}). This strategy was succsesfuly used by many authors (see Voloch, Van der Geer, Garcia, Quoos, Garzon, Borges,...).  As far as this approach is concerned, a great knowledge of the $m.v.s.p.$ can  somewhat tell us what the best possible achievement is. That gives us a motivation to further investigate such  polynomials.
%
%The properties of the polynomials will reflect directly on the geometry of the curve. So we try to establish as many prperties as we can.
%
%Any of a $m.v.s.p.$ tanking $\F_{q^n}$ into $\F_q$ gives rise to $m.v.s.p.$ tanking $\F_{q^n}$ into some subsets of $\F_q$. We characterize all such polynomials and point out some of its properties. As an aplication, we provide a generalization of the Hermitian curve whose ratio
%$N/g$ is bigger than the one in [G-Sti].

Let $q$ be a power of a prime $p$, and  for any non-constant polynomial $F\in\fq[x]$, let $V_F=\{F(\alpha): \alpha \in \F_q\}$ be its value set. Since  $F$ has no more than  $\deg F$ zeroes, one can easily show that  $V_F$ satisfies
\begin{equation}\label{mvsp-def}
\left\lfloor\dfrac{q-1}{\deg F}\right\rfloor +1 \leq |V_F| \leq q,
\end{equation}
where $\left\lfloor t \right\rfloor$ is the greatest integer $\leq t$, and $ |\mathcal{S}|$ denotes the cardinality of the set $\mathcal{S}$. If $ |V_F| $ attains the lower bound in $(\ref{mvsp-def})$, then $F$ is called a \emph{minimal value set polynomial}. Throughout this text, a  polynomial of this kind  will be referred to by the abbreviation  $\mvsp$.

%Characterizing special classes of polynomials over finite fields and determining its cardinal number  are fundamental theoretical problems which are not only worthy of study in their own right, but also because of their applications.
 
%A  particular motivation to  investigate  $\mvsps$ is related to its potential application on  the construction of new curves with many rational points.  The study connecting  some special aspects of $\mvsps$ and certain  curves with many points is being handled in a companion paper.

In 1961, motivated by a generalization of the Waring's  problem modulo $p$, the authors of \cite{carlitz_pol_minimal_61}  characterized $\mvsps$ of degree $< 2p+2$. They also remarked that for any polynomial $F\in \F_{q^{n}}[x]$, with $0<\deg F\leq (q^{n}-1)/(q-1)$ and $V_F \subset \F_q$, the polynomial $G=F^{a}$, where $a|(q-1)$, is an $\mvsp$ in  $\F_{q^{n}}[x]$ (see \cite[Section 5]{carlitz_pol_minimal_61}).
 
 In 1963,  W. Mills \cite{mills_pol_minimal_64} provided additional results on $\mvsps$ and managed to  characterize all $\mvsps$ in $\F_{q}[x]$ of degree $\leq \sqrt{q}$.  Such polynomials (see Theorem \ref{the:deg_lt_q}) turned out to be part of a larger class of $\mvsps$ previously presented in  \cite{carlitz_pol_minimal_61}. Mills also provided a complete classification of $\mvsps$ in $\F_{p^2}[x]$ and remarked that for $q>p^{2}$, an $\mvsp$ different from the ones already mentioned in \cite{carlitz_pol_minimal_61}  would be  unlikely to exist.

%They also identified the following two classes of $\mvsps$:  
%
%Let $G$ be a polynomial over $\fqn$ and let $v$ be a positive divisor of $q-1$. Then  $G^v$ is a $\mvsp$ whenever: 
%\begin{itemize}
% \item  $G$  satisfies  $0<\deg G\leq \dfrac{q^{n}-1}{q-1}$ and $V_G \subset \F_q$; or
%\item  $G$ is of the form $G=\sum\limits_{i=0}^{m}\omega_i x^{q^{i}}$ and splits over $\fqn$. 
%\end{itemize}
%
%In 1963,  W. Mills \cite{mills_pol_minimal_64} provided additional results on $\mvsps$ and managed to  characterize all $\mvsps$ in $\F_{q}[x]$ of degree $\leq \sqrt{q}$.  Such polynomials (see Theorem \ref{the:deg_lt_q}) turned out to be part of the second class of $\mvsps$ aforementioned. Mills also provided a complete classification of $\mvsps$ in $\F_{p^2}[x]$ and remarked that for $q>p^{2}$, a $\mvsp$ different from the ones above would be  unlikely to exist.

In the first part of this paper (section \ref{sec:prelim}), we extend the work of Mills in two ways. First we characterize  $\mvsps$ over $\fq$ of degree $\sqrt{q}+1$; second  we show that, up to composition with degree one polynomials, there is a unique $\mvsp$ in $\fq[x]$ of degree $\leq \sqrt{q}$ for a given set of values. %It turns out that, up to affine transformations, there is  only one such $\mvsp$ over $\fq$, namely $x^{\sqrt{q}+1}$.

 Despite the remarkable work of Mills in \cite{mills_pol_minimal_64}, which provides a great deal of information on $\mvsps$, there exists no complete solution to the following  fundamental problem: 

\begin{question}\label{question}
 Given $\mathcal{S}\subset \fq$, can we find an $\mvsp$ $F \in \fq[x]$ with  $V_F=\mathcal{S}$?  If so, how many such polynomials are there? Can they be characterized? 
\end{question}

% Since that time, not much has been done with respect to the characterization of $\mvsp$.

%As we saw in the previous section,  and it leads to the characterization of all $\mvsps$ over $\fq$ of degree $\leq \sqrt{q}+1$. 

The second part of this work is built around these questions and is organized as follows: In Section \ref{sec:mvspxqx} we will provide a complete solution to this problem in a special but very important case: when $\mathcal{S}=\ff_{q'}$ is a subfield of $\fq$. The general case is treated in Section \ref{sec:mvsp_fix_set}, where we show how the work of Mills imposes some very restrictive and necessary conditions on $\mathcal{S}$. For the sets $\mathcal{S}$ satisfying such conditions, in Section \ref{sec:vec_wts} we will be able to provide a partial answer to the question of the existence of $\mvsps$ with $\mathcal{S}$ as its set of values and give a lower bound to the number of such polynomials.  As a by-product of our constructions, we  show that $\mvsps$ other than the ones expected by Mills do exist: a simple example being the polynomial $F=x^{q^{4}+q}-x^{q^{3}+1}$ in $\F_{q^{6}}[x]$, which is addressed in Section \ref{sec:examples}. 

The characterization of special classes of polynomials over finite fields, and the determination of its cardinal number  are fundamental theoretical problems, which are not only worthy of study in their own right, but also because of their applications. Our interest in  investigating  $\mvsps$ is related to its potential application to the construction of new curves with many rational points, the characterization of certain Frobenius non-classical curves and the construction of elliptic curves over $\fq(t)$ with many integral points.  The study connecting special aspects of such polynomials and the aforementioned applications  is being handled in companion papers.

\text{}\\
{\bf Notation:} Throughout the text, $q$ will denote a power of a prime $p$.

\section{Preliminaries}\label{sec:prelim}

The content of this section is mostly based on \cite{mills_pol_minimal_64}, from where a large amount of notations and results are borrowed. Therefore, we believe that  some familiarity with Sections $1$ and $2$ of  that paper would be desirable for  our reader.
%\begin{remark}
%%%Caso $r\leq 2$ e transformacoes afins.
%It is worth mentioning that if  $F \in \F_q[x]$ is a $\mvsp$ and   $G=\alpha x+\beta \in  \F_q[x]$ is non-constant,  then $F_1=F(G(x))$ and  $F_2=G(F(x))$  are $\mvsp$ with   $V_{F_1}=V_F$ and $V_{F_2}=\alpha V_F+\beta$. This simple fact will be addressed sometimes during our  discussion.
%%\end{remark}
%%\end{enumerate}
%%\begin{enumerate}
%%\item $\tilde{F}=F(ax+b)$ is  is a $\mvsp$  with $V_{\tilde{F}}=V_F$.
%%\item $\tilde{F}=aF(x)+b$ is  is a $\mvsp$  with $V_{\tilde{F}}=aV_F+b$.
%%
%%
%\end{remark}

We begin by pointing out that $\mvsps$ $F$ with  $|V_F|\leq 2$ do not fit in the general pattern and were handled separately in \cite{carlitz_pol_minimal_61}.  Mills \cite{mills_pol_minimal_64} proved the following important result  that presents necessary conditions for a polynomial $F\in\fq[x]$ with $|V_F|>2$ to be an $\mvsp$ and encloses other valuable information about such polynomials.

\begin{theorem}[Mills]\label{main} %Let $\F_q$ be a field of characteristic $p$, and 
Let  $F \in \F_{q}[x]$ be a polynomial of positive degree. Let $V_F=\{\gamma_0,\gamma_1,\cdots,\gamma_r\}$ be its set of values, and let 
$l_i$ be the number of distinct roots in $\F_q$ of the polynomial $F-\gamma_i$. Let $\gamma_i$ be arranged in such way that $l_0\leq l_i$, $1\leq i \leq r$. Set $L=\prod (x-\pi)$ where the product is over the distinct roots of $F-\gamma_0$ that lie in $\F_q$. 

Suppose that $F$ is an $\mvsp$ and $r>1$. Then there exist positive integers $v,m,k$; a polynomial $N$ over $\F_q$, and $\omega_0,\omega_1,\cdots,\omega_r$ in $\F_q$ such that $L\nmid N$,
$v\mid (p^{k}-1)$, $1+vr=p^{mk}$, $L'$ is a $p^{mk}$-th power, $\omega_0\neq 0$, $\omega_m=1$, 
\begin{equation}\label{Fv}
F=L^{v}N^{p^{mk}}+\gamma_0,
\end{equation}

\begin{equation}\label{valueset}
\prod \limits_{i=1}^{r}(x-\gamma_i+\gamma_0)=\sum \limits_{i=0}^{m}\omega_i x^{(p^{ki}-1)/v},
\end{equation}

\begin{equation}\label{L-id}
\sum \limits_{i=0}^{m}\omega_i L^{p^{ki}}N^{p^{mk}(p^{ki}-1)/v}=-\omega_0(x^{q}-x)L'.
\end{equation}
\end{theorem}
\begin{proof}
See \cite[Theorem 1]{mills_pol_minimal_64}.
\end{proof}

As a consequence of this theorem, Mills presented the following characterization of  all   $\mvsps$ over $\fq$ of degree $\leq \sqrt{q}$. 
\begin{theorem}[Mills]\label{the:deg_lt_q} Let  $F \in \F_q[x]$ be a  polynomial with 
$0<\deg F \leq \sqrt{q}$. Then $F$ is an $\mvsp$ if and only if $F$ is of the form

\begin{equation}\label{eq:deg_sqrtq}
F=\alpha L^{v}+\gamma,
\end{equation}
where $L$ is a polynomial that factors into distinct linear factors over $\F_q$ and that has the form
\begin{equation}
L=\sum\limits_{i=0}^{d}\varphi_i x^{p^{ki}}+\beta,
\end{equation}
where $v$ and $k$ are  integers such that $v|(p^{k}-1)$, $q$ is a power of $p^{k}$, and $\alpha$, $\beta$, $\gamma$ and the $\varphi_i $ are elements of $\F_q$.
\end{theorem}
\begin{proof}
 See \cite[Theorem 2]{mills_pol_minimal_64}.
\end{proof}
In the following result, we show how Theorem \ref{main} can also be used to extend Mills' characterization to the case of degree $\sqrt{q}+1$.

%In our first  result, we use this theorem to provide a characterization of $\mvsps$ of degree $\sqrt{q}+1$.

%The next result extends Theorem \ref{the:deg_lt_q} to the case $\deg F\leq \sqrt{q}+1$.

% Our proof of this statement relies on the following fundamental theorem on $\mvsps$ proven by Mills:
\begin{theorem}\label{q+1}
Let  $F$ be a polynomial over $\fq$ with $|V_F|>2 $ and  of degree $\sqrt{q}+1$. $F$  is an $\mvsp$ if and only if there exist $\alpha,\beta,\gamma \in \F_q$, with $\alpha\neq0$, such that 
$$
F=\alpha(x+\beta)^{\sqrt{q}+1}+\gamma.$$ 
%In particular,   the value set of $(F-\gamma)/\alpha$ is $\F_{\sqrt{q}}$.
\end{theorem}
\begin{proof}
%\begin{equation}\label{Fv}
%F=L^vN^{p^{mk}}+\gamma_0
%\end{equation}
%and
%
%\begin{equation}\label{L-id}
%\sum \limits_{i=0}^{m}\omega_iL^{p^{ki}}N^{p^{mk}(p^{ki}-1)/v}=-\omega_0(x^{p^{n}}-x)L',
%\end{equation}
%
%with the definitions and assertions therein contained.

We use the notation of Theorem \ref{main}.

Suppose that $q=p^{n}$ and  that $F$ is an $\mvsp$. Since $\deg F= p^{n/2}+1$, it follows from the definition of $\mvsp$ that $|V_F|=\dfrac{p^{n}-1}{p^{n/2}+1}+1=p^{n/2}$.  Thus $r=p^{n/2}-1$ and  $1+vr=p^{mk}$ give $v=\dfrac{p^{mk}-1}{p^{n/2}-1}$. Therefore $n/2$ divides $mk$, and since $mk \leq n$ (implicitly given by \eqref{Fv} and \eqref{L-id}) we have that either $mk=n/2$ or $mk=n$.  If  $mk=n$, then $v=p^{n/2}+1$ and $(\ref{Fv})$ gives
$$F=L^{p^{n/2}+1}N^{p^{n}}+\gamma_0.$$
Since $\deg F=p^{n/2}+1$ and $\deg L \geq 1$, we have that $\alpha:=N^{p^{n}} \in \F_q^*$ and $L=x+\beta$ for some $\beta\in \F_q$. We set $\gamma:=\gamma_0$ and obtain

\begin{equation}\label{F-q+1}
F=\alpha(x+\beta)^{p^{n/2}+1}+\gamma,
\end{equation}
 as desired.

If   $mk=n/2$ then $v=1$, and we will show that such case does not occur. In fact, we suppose it does occur and use  $(\ref{Fv})$   to write

\begin{equation}\label{Fv1}
F=LN^{p^{n/2}}+\gamma_0.
\end{equation}
This implies  that either 
%and
%
%\begin{equation}\label{L-idv1}
%\sum \limits_{i=0}^{m}\omega_iL^{p^{ki}}N^{p^{n/2}(p^{ki}-1)}=-\omega_0(x^{p^{n}}-x)L',
%\end{equation}
% 
% where $L'$ is a $p^{n/2}$-th power.  Now  from $(\ref{Fv1})$ we have
 \begin{enumerate}
 \item[(i)] $L$ and $N$ are polynomials of degree $1$; or
 \item[(ii)] $\deg  L=p^{n/2}+1$  and $N$ is constant.
 \end{enumerate}

If we assume that $L$ and $N$ have degree 1, then $F-\gamma_0=LN^{p^{n/2}}$ has all of its roots defined over $\fq$. Since $L=\gcd(F-\gamma_0,x^q-x)$ has degree 1, it follows that $F-\gamma_0$ cannot have two distinct roots. Therefore $L|N$,  and this contradicts Theorem \ref{main}. 

If we assume that (ii) is true, then  $F-\gamma_0$ has $p^{n/2}+1$ distinct roots over $\ff_{p^n}$. For all $0\leq i\leq r$, recall that in Theorem \ref{main} it is assumed that $l_0\leq l_i$.  Therefore, $F-\gamma_i$ also has $p^{n/2}+1$ distinct roots over $\ff_{p^n}$. This implies that $\prod_{i=0}^r (F-\gamma_i)$ is a polynomial over $
\ff_{p^n}$ with  $(r+1)(p^{n/2}+1)=p^n+p^{n/2}$ distinct roots defined over $\ff_{p^n}$. This clear contradiction implies that neither (i) nor (ii) can hold and the result follows. The converse is trivial.
\end{proof}

 %One consequence of the previous result is that the set of values of a $\mvsp$ of degree $\sqrt{q}+1$ is of the form $\alpha\fq+\beta$, for some $\alpha\in\fq^*$ and $\beta\in\fq$. Moreover, 

Observe that this result can be rephrased as follows: up to an affine composition\footnote{\label{foot:affinecomp} It is worth mentioning that if  $F \in \F_q[x]$ is an $\mvsp$ and   $G=\alpha x+\beta \in  \F_q[x]$ is non-constant,  then $F_1=F(G(x))$ and  $F_2=G(F(x))$  are $\mvsps$ with   $V_{F_1}=V_F$ and $V_{F_2}=\alpha V_F+\beta$.  The polynomials $F_1$ and $F_2$ will be referred to as \emph{affine compositions} of $F$.}, $x^{\sqrt{q}+1}$ is the only $\mvsp$ of degree $\sqrt{q}+1$.   A similar statement  is true for  $\mvsps$ of degree $\leq\sqrt{q}$ with a given set of values: they are also unique up to affine compositions.   Our proof of this fact relies in part on  the following consequence of the results in \cite{mills_pol_minimal_64}.%, which is interesting on its own.the characterization of $\mvsps$ of degree $\leq\sqrt{q}$ provided by Mills  and 

\begin{lemma}\label{F-propert.} If  $F \in \fq[x]$ is  an $\mvsp$ with $V_F=\{\gamma_0,\gamma_1,\cdots,\gamma_r\}$ and $r>1$, then
 \begin{enumerate}
% \item [(i)] Either  $T'$ is  constant or $T=\frac{1}{k}x(T'+a_1 (k-1))$.
% \item [(ii)]If $i<k$,  then $p\nmid i$ implies $a_i=0$. 
%\item [(iii)]  $\sum\limits_{a_i \in V_F} a_i=0$; that is, $a_{k-1}=0$. 
%\item [(iv)] $T'$ is constant if and only if $\mid V_F \mid \equiv 0 \mod p$.
%\item [(i)]$T'$ is constant  if and only if $F''$ is zero. In particular, $T'$ is constant if $char=2$.
%\item [(i)] $\deg F \mid (q-1)$ if and only if $p \nmid \deg F$.
\item [(i)] A root of $F_i=F-\gamma_i=0$ is not in $\F_q$ if and only if its multiplicity is divisible by $p$.
\item [(ii)]$F-\gamma_i=0$ has a simple root for at least $r$ values $\gamma_i \in V_F$. 
\end{enumerate}
\end{lemma}
%The next result provides some facts about the polynomial $T$, along with  some aditional  properties of the polynomial $F$.
\begin{proof} 
%\begin{enumerate}
%\item [(i)] Clearly $T(F)=\theta(x^q-x)F'$ implies $T'(F)F'=\theta(x^q-x)F''-F'\theta$, and a quick inspection on the latter identity gives the result.
%\item [(i)] 
%To prove the first item, we write $T(F)=\theta(x^q-x)F'$  and obtain
%$$\deg F'=(r+1)\deg F-q= \deg F(\left[\frac{q-1}{\deg F}\right] +1) -q  \leq \deg F (\frac{q-1}{\deg F} +1) -q= \deg F-1.$$ Note that $\deg F \nmid (q-1)$ if and only if $\left[\frac{q-1}{\deg F}\right]<\frac{q-1}{\deg F}$. Thus $\deg F \nmid (q-1)$ is equivalent to $\deg F'< \deg F-1$, which in turn is equivalent to  $p \mid \deg F$.
This proof uses notations and results of \cite{mills_pol_minimal_64}. Item $(i)$ follows directly from Lemmas 2 and 3 in \cite{mills_pol_minimal_64}.  In fact, since $F_i=L_iU_i$ where $L_i=\gcd(F_i,x^q-x)$, Lemma 2 in \cite{mills_pol_minimal_64} implies $F_i=L_i^{w_i+1}H_i^p$, which proves the ``only if" part of (i). The other part of $(i)$ follows from the fact that $F_i$ is not a $p$-power; that is   $F_i'=F'\neq 0$, which can be read off 
Theorem \ref{main}.

%But $F'_i=F'\neq 0$, otherwise \eqref{Fv} would imply that either $p|v$ or $L'=0$, which contradicts Theorem \ref{main}.

Lemma 2 in \cite{mills_pol_minimal_64} also gives $L_i \nmid H_i$ which, together with Lemma 3 in \cite{mills_pol_minimal_64}, implies $(ii)$.
\end{proof}

%, the norm polynomial from $\fq$ to $\ff_{\sqrt{q}}$. 

%For the sake of completeness, let us include here Mill's classification of $\mvsps$ of degree $\leq\sqrt{q}$.

As remarked in footnote \ref{foot:affinecomp}, given $a\in\fq^*$, $b\in\fq$ and a non-constant polynomial $F$, the polynomials $F$ and  $G:=F(a x+b)$ have the same set of values.  In what follows, we will show that 
any two $\mvsps$ of degree $\leq \sqrt{q}$ and with the same value set are related in this way. %It is in this sense that we consider $\mvsps$ of degree $\leq\sqrt{q}$ with a given set of values to be unique. 

%The next result shows that . It is in this sense that we consider $\mvsps$ of degree $\leq\sqrt{q}$ and given set of values to be unique.

%In what follows, we extend this result in two different ways: first we show that, in a sense that will be explained shortly, there is a unique $\mvsp$ with degree $\leq \sqrt{q}$ and given set of values. Later, we characterize the $\mvsps$ over $\fq$ of degree $\sqrt{q}+1$.  

%Using Theorem \ref{main}, Mills  managed to  characterize the $\mvsps$  $F\in \F_q[x]$ with $\deg F \leq \sqrt{q}$. 
 %We begin by asking if, for  such a polynomial $F$, one can find another  $\mvsp$  $G\in \F_q[x]$ with  $V_G=V_F$ \textcolor{red}{Sim. Tome $G=F(ax+b)$}. 
 %As we will show shortly,  such a polynomials $F$  with  $\deg F \leq \sqrt{q}$ are somewhat unique. Before discussing  that, we state Mills' result (see \cite{mills_pol_minimal_64} ) regarding  this special class of $\mvsps$.

\begin{proposition}\label{raiz}
Let $F \in\F_q[x]$ be an $\mvsp$ such that $|V_F|>2 $ and $\deg F \leq \sqrt{q}.$ If $G \in\F_q[x]$ is an $\mvsp$ such that  $\deg G \leq \sqrt{q}$ and $V_F=V_G$  then $G=F(ax+b)$ for some  $a,b \in \F_q$.
\end{proposition}

%Given the above result, we can think of $F$ as being unique   in the following sense: For any polynomial $F \in \F_{q}[x]$, the polynomial $G=F(ax+b)$ where $a\in \F_q^{*}$ and $b\in \F_q$ satisfies
%$V_G=V_F$. That is just because $ax+b$ permutes the elements of $\F_q$.

%\noindent{\bf Proof of Lemma \ref{raiz}.}
%\textcolor{red}{Como voce mesmo observou a prova deste teorema eh mais simples e cai bem com os resultados da ultima secao. Proponho muda-lo pra secao 4 e inclui-lo na minha prova da sobrejetividade. A nao ser que voce precise deste resultado pra um outro artigo.}
\begin{proof}
Without loss of generality, we may assume $d:=\deg F \geq \deg G$.  Suppose $G(x)\neq F(ax+b) $ for all $a,b\in \F_q$; that is, $F(x)-G(y)$ has no $\F_q$-linear factors (in particular $d \geq 2$). Since $d \leq \sqrt{q}$,  Theorem \ref{the:deg_lt_q}   implies $F=\alpha L^{v}+\gamma_0$ and  so $F-\gamma=0$  has $d$ distinct roots  for all $\gamma \in \{\gamma_1,\gamma_2,\cdots, \gamma_r\}=V_F\backslash \{\gamma_0\}$, and by Lemma \ref{F-propert.}.$(i)$, all such roots lie in $\F_q$. Therefore, if $n$ is the number of times $G(\lambda)$ assumes $\gamma_0$ as $\lambda$ ranges over $\F_q$ then the number $N$ of $\F_q$-solutions to $F(y)=G(x)$ is given by $N=d(q-n)+n\cdot \deg L$. But since $1\leq n\leq d \leq \sqrt{q}$ we have
\begin{equation}\label{npontos}
N=d(q-n)+n\cdot \deg L \geq d(q-d)+\deg L\geq q(d-1)+1
\end{equation}

Notice that if $N=q(d-1)+1$, then (\ref{npontos}) gives $1=n=d$, contradicting $d\geq 2$. Therefore we have 
$$N>q(d-1)+1,$$
which contradicts the main result (Theorem 3.1) of \cite{Homma_kim_conj_plane}. This completes the proof.

\end{proof}

The results in this section show that $\mvsps$ with degree $\leq\sqrt{q}+1$ are unique (up to affine compositions) if their set of values are fixed. However, this is far from the truth in the general case. As we will see in the following sections, the set of $\mvsps$ with a fixed set  of values $\mathcal{S}$ can be quite large.
 
%\textcolor{red}{Modificar este final de secao}

%Given the relation between the polynomials in Theorem \ref{q+1} and the $\mvsp$  $F\in\F_q[x]$ with $V_F=\F_{\sqrt{q}}$, it would be of interest to give a  list the latter ones. This will be done  in more generality in the next section.

\section{On m.v.s.p.'s with a given set of values}\label{sec:mvsp_fix_set}

%\textcolor{red}{Colocar introducao para esta secao}

Our discussion of Problem \ref{question}  starts with the following criterion to decide whether or not a polynomial $F$ is an  $\mvsp$ with $V_F=\mathcal{S}$,  for a given set $\mathcal{S}\subset\fq$. The following theorem, which  can be read off  the work of \cite{carlitz_pol_minimal_61}  and \cite{mills_pol_minimal_64}, will be used to restate Problem \ref{question} in terms of the polynomial $T=\prod\limits_{\gamma \in \mathcal{S}}(x-\gamma)$ of values of $F$.

%Let us start our discussion of Problem \ref{question} by providing a characterization of $\mvsps$ with a given $\mathcal{S}\subset\fq$ as its set of values. The following criterion, which  can be read off  the work of \cite{carlitz_pol_minimal_61}  and \cite{mills_pol_minimal_64},  will be of fundamental importance to what follows. %We should mention  that  part of  the  proof  will be addressing notations and results stated in \cite{mills_pol_minimal_64}. 
\begin{theorem}\label{teoMills}
Let $F$ be a polynomial in $\F_q[x]$, and let $\mathcal{S}$ be a subset of $\F_q$. If  $|\mathcal{S}|>2$ and $T=\prod\limits_{\gamma \in \mathcal{S}}(x-\gamma)$, then the following are equivalent:
\begin{enumerate}
\item[(A)] $F$ is an $\mvsp$ with $V_F=\mathcal{S}$.
\item[(B)] There exists $\theta\in\fq^*$ such that
\begin{equation}\label{eq:minvset}
T(F)=\theta(x^q-x)F'.
\end{equation}
\end{enumerate} 
Moreover, $\theta= -T'(\gamma)$ for some $\gamma\in\mathcal{S}$.
%$\theta=-\displaystyle\prod_{j=1\atop i\neq j}^k (a_i-a_j)$, for some $1\leq i\leq k$.
\end{theorem}
\begin{proof}

If $F$ satisfies $(B)$, then we clearly have $V_F\subset \mathcal{S}$,  and thus $|V_F|\leq |\mathcal{S}|$.
Equating degrees in (\ref{eq:minvset}), we obtain  $|\mathcal{S}| \cdot\deg F\leq q+\deg F-1$, i.e.  
\begin{equation}\label{eq:degree_of_mvsp}
|\mathcal{S}|-1\leq \frac{q-1}{\deg F}
\end{equation} and thereby
$$
|\mathcal{S}|\leq \left\lfloor\frac{q-1}{\deg F}\right\rfloor+1\leq |V_F| \leq |\mathcal{S}|,
$$
which implies (A). 

That (A) implies (B) follows from Lemma 1 in \cite{mills_pol_minimal_64}.

The fact that $\theta= -T'(\gamma)$, for some $\gamma\in\mathcal{S}$, can be proved by differentiating both sides of \eqref{eq:minvset} and using the fact that $F'\neq 0$.
\qedhere% (see the definition of $\zeta$ and Lemmas 1, 2 and 3 therein). 
\end{proof}

\begin{remark}
 Using Lemmas 1, 2 and 3 in \cite{mills_pol_minimal_64}, one can actually prove that $\theta=-T'(\gamma_i)$ for all $i\neq 0$ (notation as in Theorem \ref{main}).
\end{remark}

For our convenience, let us label the following assumption 
\begin{equation}\label{H}
\textit{separable monic polynomial of degree $>2$ which splits over $\fq$}\tag{*}
\end{equation}
which will be used several times throughout the text.
%satisfied by the polynomial $T$ whose set of roots is $\mathcal{S}$

%In light of the previous result, let us define for any $T$ satisfying (\ref{H}), the set 
%\begin{equation}\label{def:wt} 
%\M(T|\fq):=\{F\in\fq[x]: T(F)=\theta (x^{q}-x)F', \text{ for some } \theta\in\fq^*\}.
%\end{equation}
%An element of this set will be called a \emph{$T-\mvsp$ over $\fq$}.

%In light of Theorem \ref{teoMills}, we make the following definition.
%The previous result motivates the  following.
\begin{definition}\label{def:wt}
Let $T$ be a  polynomial satisfying \eqref{H}. The set of \emph{$T-\mvsps$ over $\fq$} is defined to be
$$
\M(T|\fq):=\{F\in\fq[x]: T(F)=\theta (x^{q}-x)F', \text{ for some } \theta\in\fq^*\}.
$$
We will say that $\M(T|\fq)$ is \emph{non-trivial} if it contains a non-constant element.
\end{definition}

It follows from Theorem \ref{teoMills} that $F$ is a $T-\mvsp$ over $\fq$ if either $F=\alpha$ is a root of $T$ or $F$ is an $\mvsp$ whose set of values is the set of roots of $T$.  Therefore,  our investigation on $\mvsps$ with given set of values can  be done via the sets $\M(T|\fq)$ and  we can deal with the following restatement of  Problem \ref{question}: % can be rephrased as follows: 

\begin{proba}\label{prob_a}
Given a polynomial $T$ satisfying (\ref{H}), is  $\M(T|\fq)$ non-trivial? How large is $|\M(T|\fq)|$? Can we determine the whole of $\M(T|\fq)$?
\end{proba}

\begin{remark} Notice that, according to  footnote \ref{foot:affinecomp}, it suffices  to solve Problem $1a$ in the  case where $x|T(x)$.

\end{remark}

%One self-evident necessary condition for the non-triviality of the sets $\M(T|\fq)$ is the assumption \eqref{H} included in the definition of such sets. 

The following  consequence of Theorem \ref{main}  points out a  necessary condition for the non-triviality of the sets $\M(T|\fq)$,  further restricting   the class of polynomials $T$, as well as the sets $\mathcal{S}$ that can be  the value set of an $\mvsp$.
% which can be used to restrict the class of polynomials $T$  being considered.%Concerning , we should notice the following necessary condition, a consequence of :

%be limited to the cases where $T$ is given by  the right hand side of $(\ref{valueset})$.

%We will say that $\M(T|\fq)$ is \emph{non-trivial} if $\M(T|\fq)\neq  \rc(T)$, where $\rc(T)$ is defined to be the set of roots of $T$. Therefore the question of whether a set $S$ is a value set of certain $\mvsp$ is equivalent to the question of whether $\M(T_S|\fq)$ is non-trivial, for $\displaystyle T_S(x):=\prod_{\alpha\in S} (x-\alpha)$. 

%A simple example of $\M(T|\fq)$ with a non-trivial element is the following: Let $d$ be a positive integer such that $(d-1)|(q-1)$. Then it is easily verified that $x^{\frac{q-1}{d-1}}\in\M(x^d-x|\fq)$.  The case where $d=q'$ is a prime power and $\F_{q'} \subset \F_q$ is of special interest, and it will be the main object of study in the next section.

\begin{proposition}\label{pro:nec_cond}
Let $T$ be a polynomial satisfying (\ref{H}). If $\M(T|\fq)$ is non-trivial, then there exist positive integers  $k$,  $m$ and $v$, with $v|(p^k-1)$, and elements $\gamma, \omega_1,\ldots, \omega_m\in\fq$  such that 
\begin{equation}\label{eq:addpol}
A(x):=T(x^{v}+\gamma)/x^{v-1}=\sum\limits_{i=0}^{m}\omega_i x^{p^{ki}}.
\end{equation}
\end{proposition}
\begin{proof}
Let $\{\gamma_0,\ldots,\gamma_r\}$ be the set of roots of $T$. By Theorem \ref{teoMills}, any  non-constant $F\in\M(T|\fq)$ is an $\mvsp$ with $V_F=\{\gamma_0,\ldots,\gamma_r\}$ and $r>1$. Therefore, $F$ satisfies the hypotheses of Theorem \ref{main}. In particular, we can guarantee the existence of positive integers $v$, $k$ and $m$, and elements  $\omega_1,\ldots,\omega_m\in\fq$ such that  (\ref{valueset}) holds. Our conclusion is just a restatement of (\ref{valueset}) with $\gamma=\gamma_0$.
%This is a consequence of Theorem \ref{main} and equation $(\ref{valueset})$.
\end{proof}
%  $(\ref{valueset})$ points out restrictions on the  monic polynomial  $T$. In fact  notice that, up to a root-translation,  for $\M(T|\fq)$ to be non-trivial it is necessary that there exist positive integers $v$, $k$ and $m$ such that 
Polynomials given by the right-hand side of \eqref{eq:addpol}  are called \emph{$p^k$-additive} polynomials  and have been extensively studied (see \cite[Section 3.4]{lidl_finite_field_83})\footnote{In \cite{lidl_finite_field_83}, $p^k$-additive polynomials are called $p^k$-linearized polynomials.}. As we will see, the class of additive  polynomials plays a crucial role in the construction (and possibly in  the complete characterization) of $\mvsps$ with a given set of values. Most notably is the following result which is a first step towards a proof of a partial converse to the previous proposition.

\begin{proposition}\label{pro:generalmt}
Let $v$ be a positive divisor of $p^k-1$ and $T$ be  a polynomial satisfying \eqref{H}. Assume $x|T(x)$ and $A(x)=T(x^v)/x^{v-1}$ is a $p^k$-additive polynomial satisfying \eqref{H}. Then
\begin{eqnarray*}
\M(A|\fq) & \longrightarrow & \M(T|\fq)\\
F &\longmapsto & F^v
\end{eqnarray*}
is a well-defined map.
%$F\in\M(A|\fq)$ if and only if $F^v\in\M(T|\fq)$.
In particular, 
$$
\left|\M(T|\fq)\right|\geq \frac{|\M(A|\fq)|-1}{v}+1.
$$
%$\deg T> 2$
\end{proposition}
% that In the next proposition we will generalize this example and show other ways to construct non-trivial sets of $T-\mvsp$'s. The following notion will be used:
\begin{proof}
For $F\in\M(A|\fq)$, we have
$$
T(F^v)=F^{v-1}A(F)=F^{v-1}(-A'(x^{q}-x)F')=-\frac{A'}{v}(x^{q}-x)(F^v)'.
$$
Thus $F^v\in\M(T|\fq)$, by definition. %The same equality can be used to prove the converse.
\end{proof}

Propositions \ref{pro:nec_cond} and \ref{pro:generalmt} indicate that a solution to  Problem \hyperref[prob_a]{1a} in the particular case of $p^k$-additive polynomials  will  provide  a partial solution to the general problem. It turns out that whenever  $A$ is $p^k$-additive polynomial satisfying (\ref{H}), we are able to show that $\M(A|\fq)$ is non-trivial  and find a lower bound  for its cardinal number. This is accomplished in the following way:

For any $p^k$-additive polynomial $A$ satisfying (\ref{H}), it is easily verified that the set $\M(A|\fq)$ is a finite $\F_{p^k}$-vector space. In Section \ref{sec:vec_wts} we will use  this extra structure to show that for any such $A$, there exist a polynomial $x^{p^{kd}}-\alpha x\in\fq[x]$ satisfying \eqref{H} and a $\F_{p^k}$-linear map $\phi:\M(x^{p^{kd}}-\alpha x|\fq)\longrightarrow \M(A|\fq)$ with $\ker\phi\subset \fq$.  In particular, the non-triviality of $\M(A|\fq)$ and ultimately $\M(T|\fq)$,  for $T$ as in Proposition \ref{pro:generalmt}, will follow from that of $\M(x^{p^{kd}}-\alpha x|\fq)$.  The next section will be used to show that not only the sets $\M(x^{p^{kd}}-\alpha x|\fq)$ are  non-trivial, but that they can be described very explicitly.

%It turns out that the vector spaces 

 %As outlined above, this implies the non-triviality of the sets $\M(A|\fq)$ for $p^k$-additive polynomials $A$ satisfying \eqref{H}; and, ultimately  that of $\M(T|\fq)$, for $T$ as in Proposition \ref{pro:generalmt}.

%Therefore, for any $T$ satisfying the hypothesis of Proposition \ref{pro:generalmt}, we will have a  non-trivial $\M(T|\fq)$  and  a lower bound for $|\M(T|\fq)|$.

%It turns out that whenever $A$ is $p^k$-additive polynomial satisfying (\ref{H}), the set $\M(A|\fq)$ is a finite $p^k$-vector space. As it will be shown in Section \ref{sec:vec_wts}, this extra structure allows us to show that $\M(A|\fq)$ is non-trivial. Moreover, we are able to find a lower bound for the dimension of $\M(A|\fq)$.  
%
%Therefore, for any $T$ satisfying the hypothesis of Proposition \ref{pro:generalmt}, we will have a  non-trivial $\M(T|\fq)$  and  a lower bound for $|\M(T|\fq)|$.

\section{Characterization of $\M(x^{p^{kd}}-\alpha x|\fq)$}\label{sec:mvspxqx}

Let us start our study of the vector spaces $\M(x^{p^{kd}}-\alpha x|\fq)$ by noticing that since  $x^{p^{kd}}-\alpha x$ satisfies \eqref{H}, it follows that $q=(p^{kd})^n$ for some integer $n\geq 1$ and that $\M(x^{p^{kd}}-\alpha x|\fq)$ is isomorphic to $\M(x^{p^{kd}}-x|\fq)$. Therefore, it is sufficient to consider the set $\M(x^q-x|\fqn)$. Observe that Theorem \ref{teoMills} implies that
%case 
\begin{equation}\label{M(xq-x)}
\M(x^q-x|\fqn)=\{F\in\fqn[x]: F^{q}-F=(x^{q^{n}}-x)F'\},
\end{equation}
since $\theta=-T'(\gamma)=1$, for some $\gamma\in\fq$.

%Recall that Theorem \ref{teoMills} guarantees that the set $\M(x^q-x|\fqn)$ contains all $\mvsps$ with value set $\fq$ only when $q>2$. The next lemma, among other things, will also show that the same is true for $q=2$.
Recall that $\deg T >2$ was required in the definition of  $\M(T|\fq)$, and that was motivated by the condition $r>2$ in Theorem $\ref{teoMills}$.  However, we want to stress that $q=2$ will be allowed in $(\ref{M(xq-x)})$, and the reason why this can be done is just one of the facts given by the following Lemma.

\begin{lemma}\label{lemaq=2}
If $F \in \fqn[x]$ is  non-constant with $V_F \subseteq \fq$, then the following are equivalent:
\begin{enumerate}
\item[(1)]$F$ is an $\mvsp$ with  $V_F=\fq$.
\item[(2)] $F^{q}-F=(x^{q^{n}}-x)F'$.
\item[(3)]$q^{n-1}\leq \deg F \leq (q^{n}-1)/(q-1).$
\item[(4)] $0< \deg F \leq (q^{n}-1)/(q-1).$
\end{enumerate}
\end{lemma}
\begin{proof}
\begin{itemize}
\item {\it (1) ${ \Longrightarrow}$ (2)}: In the notation of Theorem \ref{teoMills} we have $T=x^{q}-x$, and the result will follow  if $q>2$. Assuming $q=2$, we have 
\begin{equation}\label{caseq=2}
F^{2}-F=(x^{2^{n}}-x)G
\end{equation}
 for some $G\in \F_{2^{n}}[x]$, and we claim that $G=F'$. In fact, differentiating both sides of $(\ref{caseq=2})$  gives $F'=G+(x^{2^{n}}-x)G'$, and the result will follow if  $G'=0$.  By hypothesis, we have $\lfloor (2^{n}-1)/\deg F\rfloor +1=|V_F|=2$ which  implies
 \begin{equation}\label{hipoq=2}
\deg F \leq 2^{n}-1,
\end{equation}
  and thus  $(\ref{caseq=2})$ gives $\deg G\leq 2^{n}-2$; that is, $\deg xG\leq 2^{n}-1$.   So if $G'\neq 0$, then a monomial of $G$ of  odd degree, say $\alpha x^{2t+1}$, will give rise to a monomial $\alpha x^{2^{n}+2t+1}$ on the right side of $(\ref{caseq=2})$. But this will imply $\deg F > 2^{n}$, contradicting  $(\ref{hipoq=2})$. Therefore $G=F'$.
\item {\it (2) ${ \Longrightarrow}$ (3)}: Since $ F^{q}-F=(x^{q^{n}}-x)F'$ gives $q\cdot\deg F=q^{n}+\deg F'$, the result will follow from   $0\leq \deg F'\leq \deg F -1$.
\item {\it (3) ${ \Longrightarrow}$ (4)}: obvious.
\item {\it (4) ${ \Longrightarrow}$ (1)}$: 0<\deg F \leq (q^{n}-1)/(q-1)$  gives $(q^{n}-1)/(\deg F) \geq q-1$, and thus
$\lfloor (q^{n}-1)/\deg F\rfloor +1\geq q\geq |V_F|$. But since $|V_F| \geq \lfloor (q^{n}-1)/\deg F\rfloor +1$ always holds we have $|V_F| =\lfloor (q^{n}-1)/\deg F\rfloor +1$.\qedhere
\end{itemize}
\end{proof}

%Since $q={q'}^n$ for some integer $n\geq 1$, 
Now we proceed to  characterize the polynomials $F \in \M(x^q-x|\fqn)$. We begin with the larger set of polynomials $F$ over $\fqn$ satisfying
\begin{equation}\label{eq:modxqn}
F^{q}\equiv F\mod (x^{q^n}-x) \text{ and } \deg F\leq q^n-1.
\end{equation}

%Let $F$ be a polynomial satisfying \eqref{eq:modxqn}. 
If we write $F= m_1+m_2+\cdots+m_l,$  where $m_i$ are monomials of 
 $F$ with distinct degrees, then  
$$
F= (m_1^q\mod (x^{q^n}-x))+(m_2^q\mod (x^{q^n}-x))+\cdots+(m_l^q\mod (x^{q^n}-x)) .
$$
 
In other words, if $M=\{m_1, m_2,\cdots,m_l\}$, then  the map $\sigma : M \longrightarrow M$ given by $\sigma(m_i)=m_i^q\mod (x^{q^n}-x)$  is a permutation on $M$. Moreover, the following holds:

\begin{proposition}\label{pro:gaction} Let $F$ be a polynomial satisfying \eqref{eq:modxqn} and let $M$ be the set of all monomials of $F$ with distinct degrees. If $g$ is a generator of $\gal(\F_{q^n}|\F_q)$
then  
$$
g^{j}(m):= m^{q^{j}}\mod (x^{q^n}-x),
$$
defines an  action of $\gal(\F_{q^n}|\F_q)$ on $M$.
\end{proposition}
\begin{proof} 
It is a routine check and  we leave it to the reader.
%Since $f=f^q \mod (x^{q^n}-x)$, we know from the previous Remark that $m^q \mod (x^{q^n}-x) \in M$, and therefore
%$m^{q^{i}} \mod (x^{q^n}-x) \in M$ for $i=1,\cdots,n$.
%
%
%\begin{enumerate}
%\item[(i)] $a\mapsto a^{q^{n}}$ is the indentity of $G$. For a monomial $m=\alpha x^d$ of $f(x)$, we have
%$m^{q^n} \mod (x^{q^n}-x)= \alpha (x^d)^{q^n} \mod (x^{q^n}-x)=\alpha (x^{q^n})^{d} \mod (x^{q^n}-x)=\alpha x^{d} \mod (x^{q^n}-x) =m \mod (x^{q^n}-x).$
%
%\item[(ii)] For $g_1,g_2 \in G$, and a monomial $m$ of $f(x)$, it is also easy to check that 
%$(g_1\circ g_2)(m)=g_1(g_2(m))$.
%\end{enumerate}
\end{proof} 

The following  simple but useful fact is a consequence of the previous result:

\begin{lemma}\label{fatos} Let $F$ be a polynomial satisfying \eqref{eq:modxqn}. Then:
\begin{enumerate}
\item [(a)] If $m^q\mod (x^{q^n}-x)$ is a monomial of $F$, then so is $m\mod (x^{q^n}-x)$.
\item [(b)] If $m=\alpha x^{kq+1}$ is a monomial of $F$, then so is $\tilde{m}=\alpha^{1/q}x^{k+q^{n-1}}$.
%\item [(c)] If  $m$ is a monomial of $F$ and $\deg(m)< q^{n-1}$, then  $m^q$ is a monomial of $F$ distinct from $m$.
\end{enumerate}
\end{lemma}
\begin{proof}
Both cases follow directly from that fact that if $m$ is a monomial of $F$, then so is $m^{q^{n-1}} \mod (x^{q^{n}}-x)$.
\end{proof}

%Since there is a unique monomial of $F$ of a given degree, we may interpret the above action on monomials as an action of $\gal(\fqn|\fq)$ on the exponents of non-zero monomials of $F$. 
Let  $O_m$ denote the orbit of a monomial $m=\alpha x^k$ of $F$ by the action of $\gal(\fqn|\fq)$.  Let $s(m)$ denote the size of $O_m$. Then
$$
(\alpha x^{k})^{q^{s(m)}}\mod (x^{q^n}-x) =\alpha x^k,
$$
implies  that for any $m_0\in O_m$, $m_0$ is defined over $\ff_{q^{s(m)}}$. This observation is the last ingredient  in the characterization of the polynomials $F$ over $\fqn$ satisfying $F^q\equiv F\mod (x^{q^n}-x)$ and $\deg F\leq q^n-1$.
\begin{proposition}\label{propparcial} Every $F$ satisfying \eqref{eq:modxqn} can be written as $F=\sum\limits_{i=1}^{e} F_i$, where \begin{equation}
F_i=\sum_{m \in O_{m_i}} m
\end{equation}
and $m_1,\ldots,m_e$ are certain monomials of $F$. Moreover, for all $1 \leq i \leq e$, $F_i$ is defined over $\ff_{q^{s(m_i)}}$.
\end{proposition}

\begin{proof} Pick a monomial $m_1$ of $F$. The orbit $O_{m_1}$
is a set of distinct monomials of $F$.  Let $F_1$ be the sum of all such monomials. If those are all the monomials of $F$, then we are done; otherwise we can pick a monomial $m_2$ of $F$ not lying on the set $O_{m_1}$, and consider the polynomial $F_2$, sum of the monomials in the orbit $O_{m_2}$. We keep repeating this process, and it is clear that we will finish it after a finite number of steps, let us say $e$ steps. Now $F=\sum\limits_{i=1}^{e} F_i$ follows from the construction and from the fact that the sets $O_{m_i}$'s are mutually disjoint.
%\footnote{Notice that the $M_i$'s are the orbits of the underlying action of $\gal(\F_{q^n}|\F_q)$ on the set of monomials of $F$.}. $F$ satisfies $
\end{proof}

\begin{remark}\label{rem:fq_valued_pol}
Notice that a polynomial $F$ satisfies $F^q\equiv F\mod (x^{q^n}-x)$ if,  and only if,  $V_F\subset \fq$. Therefore the previous result characterizes all polynomials $F$ over $\fqn$  with $V_F\subset \fq$ and $\deg F\leq q^n -1$.
\end{remark}

The complete description of  $F\in \M(x^{q}- x|\fqn)$ will follow once we know what the monomials of such an $F$ look like. The additional information, namely $F^q-F =(x^{q^n}-x)F'$, will answer this question.

\begin{lemma} For all $F \in\M(x^{q}- x|\fqn)$,  there exist polynomials $A$ and $B$ such that  $F=A+B$ and 
\begin{equation}\label{split1}
F=x\left(\dfrac{A}{x^{q^{n-1}}}\right)^q+B^q
\end{equation}
\end{lemma}
\begin{proof}  
Let $A,B \in \F_{q^n}[x]$ be such that $F=A+B$ and  $A$ is formed by all the monomials of $F$ of degree at least $q^{n-1}$. Note that  Lemma \ref{lemaq=2}$.(3)$ implies $A \neq 0$. Then $A^q+B^q-F=x^{q^n}F'-xF'$, and after comparing degrees we have 
\begin{enumerate}
\item [(a)]  $ A^q=x^{q^n}F'$, which gives $F'= \left(\dfrac{A}{x^{q^{n-1}}}\right)^q$,
and 
\item [(b)] $ (B^q-F)=-xF'$, which gives $F=xF'+ B^q$.
\end{enumerate}

Combining $(a)$ and $(b)$, we obtain the desired result.
\end{proof}

\begin{lemma}\label{lem:power} If  $\alpha x^k$ is a non-zero monomial  of $F \in \M(x^{q}- x|\fqn)$, then 
$$
k=a_{n-1}q^{n-1}+a_{n-2}q^{n-2}+\cdots+a_1q+a_0
$$
where $a_i\in\{0,1\}$. 
%In particular, $F$ has at most $2^n$ non-zero terms. 
\end{lemma}
\begin{proof} 
We can certainly write $k=\sum\limits_{i=1}^n a_i q^{i}$, where $a_i \in \{0,1,\cdots,q-1\}$. From (\ref{split1}) we clearly have $a_0\in \{0,1\}$. If $a_0=0$ then Lemma \ref{fatos} implies that $\alpha^{q^{n-1}}x^{\frac{k}{q}}$ is a
 term of $F$ and then we  fall into an equivalent problem with fewer variables $a_i$. Therefore, we may assume $a_0=1$ and again Lemma \ref{fatos} implies that  $\tilde{\alpha} x^{q^{n-1}+\frac{k-1}{q}}$ is a term of  $F$. The new term $\tilde{\alpha}x^{q^{n-1}+a_{n-1}q^{n-2}+\cdots+a_2q+a_1}$ gives $a_{1}\in \{0,1\}$ and we repeat the process to reach $a_{2}$. After at most $n$ iterations,  we find that all $a_i$ are either zero or one.
\end{proof}

Before finishing the complete description of $\M(x^q-x|\fqn)$, let us recall that $s(m)$ denotes the size of the orbit of the monomial $m$ by the action of $\gal(\fqn|\fq)$ as described in Proposition \ref{pro:gaction}.

\begin{theorem}\label{the:charmxq}
A polynomial $F$ lies in  $\M(x^{q}- x|\fqn)$ if, and only if,  it  can be written as a sum of polynomials of the form
\begin{equation}\label{eq:char_mvsp_fq}
\sum_{i=0}^{s(m)-1} [m^{q^i} \mod (x^{q^n}-x)],
\end{equation}
where $m\in\ff_{q^{s(m)}}[x]$ is a monomial and $\deg m=a_{n-1}q^{n-1}+a_{n-2}q^{n-2}+\cdots+a_{1}q+a_0$ with $a_i\in\{0,1\}$.
%\begin{equation}\label{eq:char_mvsp_fq}
%\sum_{j=0}^{s(k)-1} \left[\alpha^{q^j} x^{kq^j}\mod (x^{q^n}-x)\right],
%\end{equation}
%where $\alpha\in\ff_{q^{s(k)}}$ and $k=k_1q^{n-1}+k_2q^{n-2}+\cdots+k_{n-1}q+k_n$ with $k_i\in\{0,1\}$.
\end{theorem}
\begin{proof}
That  such an $F\in \M(x^{q}- x|\fqn)$ can be written as sum of polynomials given by \eqref{eq:char_mvsp_fq} is a consequence of Proposition \ref{propparcial} and Lemma \ref{lem:power}. To prove the converse, we can use the fact that  $\M(x^q-x|\fqn)$ is an $\fq$-vector space and just prove that the polynomials given by \eqref{eq:char_mvsp_fq} lie in  $\M(x^q-x|\fqn)$.

Let $H$ be a polynomial given by \eqref{eq:char_mvsp_fq}. If $\deg H=0$ then $H\in\fq$. Thus we may assume that $\deg H>0$. We know that $H$ satisfies \eqref{eq:modxqn} and thereby $V_H\subset \fq$. Moreover, $0<\deg H\leq 1+q+q^2+\ldots +q^{n-1}=\frac{q^n-1}{q-1}$ and Lemma \ref{lemaq=2}$.(4)$  implies $H\in\M(x^q-x|\fqn)$.
\end{proof}

\subsection{\bf The dimension of $\M(x^q-x|\fqn)$}\label{sec:dimxqx}

Let $g$ be a generator of $\gal(\fqn|\fq)$ and $(a_0,\ldots,a_{n-1})\in\{0,1\}^n$. Then the map given by  $g(a_0,\ldots,a_{n-2},a_{n-1})=(a_{n-1},a_0,\ldots,a_{n-2})$ defines an action of $\gal(\fqn|\fq)$ on the set $\{0,1\}^n$. Let us identify the $n$-tuple $(a_0,\ldots,a_{n-2},a_{n-1})$ with the $q$-ary expansion of an integer $k=a_{n-1}q^{n-1}+\cdots+a_{1}q+a_0$, and choose representatives $k_1,\ldots,k_b$ on each orbit of the set $\{0,1\}^n$. For each $i\in\{1,\ldots,b\}$, let $d_i$ be the size of the orbit  $O_{k_i}$ and 
\begin{eqnarray*}
\V_i:=\left\{\sum_{j=0}^{d_i-1} (\alpha x^{k_i})^{q^j}\mod (x^{q^n}-x): \alpha \in\ff_{q^{d_i}}\right\}.
\end{eqnarray*}

\begin{theorem}\label{the:dimxqx}
For $i\in \{1,\cdots,b\}$, the set  $\V_i$ is an  $\fq$-vector subspace of $\M(x^q-x|\fqn)$ of dimension $d_i$. Moreover,
$$
\M(x^q-x|\fqn)=\bigoplus_{i=1}^b \mathcal{V}_i.
$$
and
$$
\dim_{\fq}\M(x^q-x|\fqn)=2^n.
$$
\end{theorem}
\begin{proof}

Theorem \ref{the:charmxq} implies that each $\V_i\subset\M(x^q-x|\fqn)$. Also, it can be easily verified that $\V_i$ is an $\fq$-vector space and that $|V_i|=q^{d_i}$. Consequently $\dim_{\fq} \V_i=d_i$.

Notice that if  $k=a_{n-1}q^{n-1}+a_{n-2}q^{n-2}+\cdots+a_{1}q+a_0$ is identified with an $n$-tuple $(a_0,a_1,\ldots,a_{n-1})\in\{0,1\}^n$, then the degree of $(\alpha x^k)^q \mod (x^{q^n}-x)$ is equal to $a_nq^{n-1}+a_1q^{n-2}+\ldots +a_{n-1}$; and, under the above identification, it is given by  the cyclic permutation $(a_{n-1},a_0,a_1,\ldots,a_{n-2})$. 

Therefore the set of exponents of monomials of elements in $\V_i$ is equal to the orbit of $k_i$ (identified with an $n$-tuple) by the action of $\gal(\fqn|\fq)$ on $\{0,1\}^n$. This implies that all elements of $\V_i$ will not share any monomial with any element in a distinct $\V_j$. Thus for all $i$,
$$
\V_i\cap (\V_1+\ldots+\V_{i-1}+\V_i+\ldots +\V_b)={0}.
$$
Theorem \ref{the:charmxq} implies that every element of $\M(x^q-x|\fqn)$ can be written as sum of elements in distinct $\V_i$'s. Therefore, $\M(x^q-x|\fqn)=\displaystyle\bigoplus_{i=1}^b \V_i$.

Now, if we let $o(d)$ denote the number of $\V_i$'s of dimension $d$, then 
$$
\dim_{\fq}\M(x^q-x|\fqn)=\sum_{d|n} d\cdot o(d).
$$
Notice that by construction $o(d)$ also counts the number of orbits of $\{0,1\}^n$ of size $d$. Thus, $d\cdot o(d)$ is equal to the cardinality of the union of all orbits of $\{0,1\}^n$ with size $d$. Since $\{0,1\}^n$ can be partitioned into orbits of size dividing $n$, we conclude that $\displaystyle\sum_{d|n} d\cdot o(d)=2^n$.
\end{proof}

\begin{corollary}\label{pro:dimfqd} %be such that $x^{q^d}-\alpha x$ splits completely over $\fqn$. 
If $x^{q^d}-\alpha x$ satisfies \eqref{H} then 
$$
\dim_{\ff_{q}}\M(x^{q^d}-\alpha x|\fqn)= d2^{n/d}.
$$
\end{corollary}
\begin{proof}
We leave the verification to the reader.
%Since $x^{q^d}-\alpha x$ splits  over $\fqn$, we have  that $\beta^{q^d-1}=\alpha$, for some $\beta\in\fqn^*$. 
%This yields an $\fq$-isomorphism $\M(x^{q^d}-x|\fqn)\longrightarrow \M(x^{q^d}-\alpha x|\fqn)$ given by $F\longmapsto \beta F$.
%%Therefore, there exists an $\fq$-isomorphism $\M(x^{q^d}-x|\fqn)\longrightarrow \M(x^{q^d}-\alpha x|\fqn)$ defined by $F\longmapsto \beta F$.
%Also, since  
%$$
%\dim_{\fq} \M(x^{q^d}-x|\fqn)=d \cdot\dim_{\fqd} \M(x^{q^d}-x|\fqn)=d \cdot\dim_{\fqd} \M(x^{q^d}-x|\ff_{{(q^d)}^{n/d}}),
%$$
%The result follows  from Theorem \ref{the:dimxqx}.
\end{proof}
\section{On the vector spaces $\M(A|\fq)$}\label{sec:vec_wts}

Let $A$ be a $p^k$-additive polynomial satisfying (\ref{H}). The aim of this section is to provide a method to construct non-constant elements of $\M(A|\fq)$, compute a lower bound for its dimension,  and point out cases where  this  bound is  sharp. As before, it is enough to consider  the $\fq$-vector spaces $\M(A|\fqn)$, where $A$ is a $q$-additive polynomial.  %The objective of this section is to provide a method to construct non-constant elements of $\M(T|\fqn)$, compute a lower bound for its dimension,  and point out cases where  this  bound is  sharp. %The first special case is given by the following  consequence of Theorem \ref{the:dimxqx}.
%is an $\fq$-vector space whose base can be easily computed whenever $q^n$ is not too large. Indeed, if we let  $a:=(q^n-1)/(\deg T-1)$ and  to be the linear map from the $\fq$-vector space of polynomials of degree $\leq a$ to the $\fq$-vector space of polynomials of degree $\leq a\cdot\deg T$, then $\M(T|\fqn)=\ker \Gamma$ can be computed via matrices.

% $T$ $q$-additive, to construct a linear map
%The results of the previous section indicate that our understanding of the vector spaces $\M(x^q-\alpha x|\fqn)$ is complete. The next couple of lemmas allow us to construct linear maps from such vector spaces to $\M(T|\fqn)$. As we will see, in many cases these linear maps will completely describe the set $\M(T|\fqn)$. Recall that $\rc(T)$ denotes the set of roots of $T$ and that it will be an $\fq$-vector subspace of $\fqn$ of dimension $t$, whenever $T$ is a $q$-additive polynomial of $\deg T=q^t$.

 The following two lemmas will be used  to construct $\fq$-linear maps $\phi: \M(x^{q^d}-\alpha x|\fqn) \longrightarrow \M(A|\fqn)$.  As discussed briefly  in the end of Section \ref{sec:mvsp_fix_set},  we will take advantage of our understanding of $\M(x^{q^d}-\alpha x|\fqn)$ to give information on $\M(A|\fqn)$. %As a matter of fact, there will be cases where   these linear maps  completely characterize the sets $\M(A|\fqn)$. 
 
Let $\rc(A)$ denote the set of roots of a polynomial $A$. Then  $\rc(A)$  is an $\fq$-vector subspace of $\fqn$ of dimension $t$, whenever $A$ is a $q$-additive polynomial satisfying \eqref{H} with $\deg A=q^t$.

\begin{lemma}\label{lem:linmaps}
Let $A,L$ and $M$ be $q$-additive polynomials. Suppose that $A$ and $L$  satisfy \eqref{H} and that $L=\gamma A(M(x))$, for some $\gamma\in\fqn^*$. % and that $T$ and $L$ splits completely over $\fqn$.  
 Then the map $\phi_M:\M(L|\fqn)\longrightarrow \M(A|\fqn)$ given by $\phi_M(F):=M(F)$ is  $\fq$-linear  with kernel $\rc(M)$.
\end{lemma}
\begin{proof}
The assertions are simple and can be easily checked. %To help fixing notation, we show  $F\in\M(T|\fq) \Longrightarrow M(F)\in\M(L|\fq)$.  In fact, if $L'=\alpha$ and $M'=\beta$, we have $T'=\alpha\beta\gamma$ and then 
%It should be clear that the map is $\fq$-linear and has kernel $\rc(M)$. We only need to check that $M(F)\in\M(L|\fq)$ whenever $F\in\M(T|\fq)$. This follows from the identity
%$$
%L(M(F))=T(F)=-T'(x^q-x)F'=-L'M'(x^q-x)F'=-L'(x^q-x)(M(F))'.
%$$
%$$
%\gamma L(M(F))=T(F)=-\alpha\beta\gamma(x^q-x)F'=-\gamma\alpha(x^q-x)(M(F))'
%$$
%gives the result.
\end{proof}

\begin{lemma}\label{the:composition}
Let  $A$ and $x^{q^d}-\alpha x$ be $q$-additive polynomials satisfying \eqref{H}. %with $\deg T=q^t$. $\alpha\in\fqn^*$ be such that $x^{q^d}-\alpha x$ splits completely over $\fqn$ and . 
If $A$ divides $x^{q^d}-\alpha x$ then there exists a $q$-additive polynomial $M$ splitting over $\fqd$ and an element $\gamma$ of $\fqn^*$ such that
%$$
%x^{q^d}-\alpha x=T(A(x))=A(T(x)).
%$$
$$
x^{q^d}-\alpha x=\gamma A(M(x)).
$$
\end{lemma}
\begin{proof}
 
Suppose $\deg A=q^t$, and assume $\alpha=1$. By assumption, the set $\mathcal{R}(A)$  is an $\fq$-vector subspace of codimension $d-t$ in $\F_{q^d}$. Therefore, from Corollary $2.5$ of \cite{Garcia_Ozbudak_2007}, there exists a unique monic $q$-additive polynomial $M \in \F_{q^d}[x]$ splitting over $\fqd$ and of degree $q^{d-t}$, such that $\mathcal{R}(A)=\{M(y): y \in \F_{q^n}\}$. But this implies that $A(M(x))\in \F_{q^d}[x]$ is monic of  degree $q^d$ and vanishes in $\F_{q^d}$. That is, 
$$A(M(x))=x^{q^d}-x.$$

%If 
%
%
%If $T$ divides $x^{q^d}-x$ then $\rc(T)$ is an  $\fq$-vector subspace of $\fqd$. Write $\fqd=\rc(T)\oplus \mathcal{U}$, for some $\fq$-vector subspace $\mathcal{U}$. Notice that if $\deg T=q^t$ then $\dim_{\fq} \mathcal{U}= d-t$. 
%
%There exists a $q$-additive polynomial $A_0$ defined over $\fqd$ such that $\rc(A_0)=\mathcal{U}$, see \cite[Theorem 3.52]{lidl_finite_field_83}. If we consider $A_0$ as a linear map $A_0:\fqd\longrightarrow \fqd$, then $\ker A_0=\rc(A_0)$ and $A_0(\fqd)=\rc(T)$. This implies that $T(A_0(\alpha))=0$ for all $\alpha\in\fqd$. In particular, we have that $x^{q^d}-x$ divides $T(A_0)$. 
%
%On the other hand, $\deg T(A_0)=q^d$, since $\deg A_0=q^{d-t}$. Thus, we can find an $\alpha\in\fqd^*$ such that $\alpha T(A_0(x))=x^{q^d}-x$. If we let $A=A_0(\alpha x)$, then $T(A)=x^{q^d}-x$. By replacing $T$ with $A$ in the above discussion, it is not hard to see that we will have $x^{q^d}-x=A(T(x))$, as desired.

For  $\alpha\neq 1$, we consider $\beta\in\fq^*$ such that $\alpha=\beta^{q^d-1}$. Since $A(x)$ divides $x^{q^d}-\alpha x$, we have that $\beta^{-q^t}A(\beta x)$ is monic and divides $x^{q^d}-x$. As proved above, there exists a unique monic $q$-additive polynomial $L$ such that $x^{q^d}-x=\beta^{-q^t}A(\beta L(x))$. If we take $M(x)=\beta L(\beta^{-1}x)$ and $\gamma=\beta^{q^d-q^t}$, then the result follows.
\end{proof}

We are ready to derive the main result of this section.
\begin{theorem}\label{the:dimMT}
%Let $d>0$ be an integer. 
Suppose $A$ and $x^{q^d}-\alpha x$ are $q$-additive polynomials satisfying \eqref{H} and that $\deg A=q^t$.  

If $A$ divides $x^{q^d}-\alpha x$ then there exists a $q$-additive polynomial $M$ with $\deg M=q^{d-t}$ such that the sequence
\begin{equation}\label{eq:exact_seq}
0\longrightarrow \rc(M) \xrightarrow{\iota} \M(x^{q^d}-\alpha x|\fqn)  \xrightarrow{\phi_M} \M(A|\fqn)
\end{equation}
is  exact,  where $\iota$ is the inclusion map and $\phi_M$ is defined by $\phi_M(F)=M(F)$. Moreover,
$$
\dim_{\fq} \M(A|\fqn)\geq d2^{n/d}-d+t,
$$
with equality if  and only if  $\phi_M$ is surjective.
\end{theorem}
\begin{proof}
Since $A$ divides $x^{q^d}-\alpha x$, Lemma \ref{the:composition} guarantees the existence of a $q$-additive polynomial $M$ satisfying 
$$
x^{q^d}-\alpha x=\gamma A(M(x)),
$$
for some $\gamma\in\fqd^*$. Thus, Lemma \ref{lem:linmaps} gives us the desired exact sequence.

The statement about the dimension follows from the Rank-nullity theorem and Corollary \ref{pro:dimfqd}.
\end{proof}

Let us recall that our study of the vector spaces $\M(A|\fqn)$ was motivated by Problem \ref{question}. The previous result, footnote \ref{foot:affinecomp} and Theorem \ref{teoMills}
show  that if $\mathcal{S}$ is any affine transformation of an $\fq$-vector space $\mathcal{V}$, then we have a positive solution to part of Problem 1: there are $\mvsps$ $F$ over $\fqn$ such that $V_F=\mathcal{S}$. Moreover, these results provide the lower bound $q^{d2^{n/d}-d+t}$ for the number of such polynomials, where  $d>0$ is an integer depending on the monic $q$-additive polynomial $A$ of degree $t>2$ such that $\mathcal{V}=\rc(A)$. In many cases this lower bound can be attained if, in Theorem \ref{the:dimMT}, we choose the unique multiple of $A$ of the form $x^{q^d}-\alpha x$ satisfying \eqref{H} and with minimal degree. When $A$ itself is of the form $x^{q^d}-\alpha x$, this is a consequence of Corollary \ref{pro:dimfqd}.  The next result shows that it also holds if $\dim_{\fq}\mathcal{V}\geq n/2$.

\begin{theorem}
%Let $d>0$ be an integer and $\alpha\in\fqn$. 
Suppose $A$ and $x^{q^d}-\alpha x$ are $q$-additive polynomials satisfying \eqref{H} and that $\deg A=q^t\geq q^{n/2}$.

If among all multiples of $A$ of the form $x^{q^d}-\alpha x$ we choose the one with minimal degree, then the sequence given by Theorem \ref{the:dimMT} can be extended to the exact sequence
$$
0\longrightarrow \rc(M) \xrightarrow{\iota} \M(x^{q^d}-\alpha x|\fqn)  \xrightarrow{\phi_M} \M(A|\fqn) \longrightarrow 0.
$$
In particular,
$$
\dim_{\fq} \M(A|\fqn) = d2^{n/d}-d+t.
$$
\end{theorem}
\begin{proof}

Let $G\in\M(A|\fqn)$. From (\ref{eq:degree_of_mvsp}) and $\deg A\geq q^{n/2}$, it follows that $\deg G\leq \frac{q^n-1}{\deg A-1}\leq  q^{n/2}+1$.

First assume that $\deg G=q^{n/2}+1$ for some $G\in\M(A|\fqn)$. Theorem \ref{q+1} implies that there exist $\alpha,\beta,\gamma\in\fqn$ such that $G=\alpha(x+\beta)^{q^{n/2}+1}+\gamma$. Since $\rc(A)=V_G=\alpha \ff_{q^{n/2}}+\gamma$ is a  vector space, we can conclude that $\gamma=0$ and $A=x^{q^{n/2}}-\alpha^{q^{n/2}-1} x$. As observed above, for $q$-additive polynomials of this form the result is a consequence of Proposition \ref{pro:dimfqd}.

Now suppose that $\deg G<q^{n/2}+1$ for all $G\in\M(A|\fqn)$. The fact that  $\M(x^{q^{n/2}}-\alpha x|\fqn)$ is $\fq$-isomorphic to $\M(x^{q^{n/2}}- x|\fqn)$ and Theorem \ref{the:charmxq} imply that $A$ is not of the form $x^{q^{n/2}}-\alpha x$.  Therefore $x^{q^n}-x$ is the multiple of $A$ of such form and with minimal degree. Theorem \ref{the:dimMT}  guarantees the existence of a $q$-linear polynomial $M$ and an exact sequence %given by \eqref{eq:exact_seq}
$$
0\longrightarrow \rc(M) \longrightarrow \M(x^{q^n}- x|\fqn)  \xrightarrow{\phi_M} \M(A|\fqn).
$$
We would like to show that the map $\phi_M(F):=M(F)$ is surjective. Since $\M(x^{q^n}-x|\fqn)=\{\alpha x +\beta: \alpha,\beta\in\fqn\}$, it follows that $\phi_M$ is surjective if $G\in\M(A|\fqn)$ implies that $G(x)=M(\alpha x +\beta)$ for some $\alpha,\beta \in\fqn$. This last assertion follows from Lemma \ref{raiz}.\qedhere

\end{proof}

%We will use the next section, to provide more detailed examples.
\section{Some Examples}\label{sec:examples}

As described in the introduction, the authors of \cite{carlitz_pol_minimal_61} identified two classes of $\mvsps$ over $\fqn$. Let $G$ be one of the following polynomials over $\fqn$: 
\begin{itemize}
 \item  a polynomial satisfying  $0<\deg G\leq \dfrac{q^{n}-1}{q-1}$ and $V_G \subset \F_q$; or
\item  a $q$-additive polynomial that splits over $\fqn$. 
\end{itemize}
If  $v$ is a positive divisor of $q-1$, $\alpha\in\fq^*$ and $\beta,\gamma\in\fq$, then  $\alpha G(x+\beta)^v+\gamma$ is an $\mvsp$. 

At the end of \cite{mills_pol_minimal_64}, the author says that ``it seems unlikely that there are any more" $\mvsps$ other than the ones  above. In this section we   show how the previous results  can be used to construct many examples of $\mvsps$; some of which are not of the form given above.

%We use the next examples to illustrate the results of the previous sections.

%\textcolor{red}{Esse seria o talvez esqueleto da secao. Agora precisa completar/ajustar os detalhes com a notacao que vc criou.}

\begin{enumerate}
\item[(1)] Obtaining the elements of $\M(x^{q}-x|\fqn)$ for  $n=3$. 

\bigskip
%\begin{itemize}
%5\item The case $n=2$.

%According to Theorem \ref{the:charmxq}, the set of exponents for the monomials of $F\in \M(x^{q}-x|\F_{q^{3}})$ is given by $\{0,1,q,q+1,q^2,q^2+1,q^2+q,q^2+q+1\}.$
Let us consider the action of $\gal(\ff_{q^3}|\fq)$ on the set of $3$-uples $(a_0,a_1,a_2)\in\{0,1\}^3$ via cyclic permutation. The orbits of this action are given by
\begin{itemize}
\item $O_0:=\{(0,0,0)\}$,
\item $O_1:=\{(1,0,0),(0,1,0),(0,0,1)\}$,
\item $O_{q+1}:=\{(1,1,0),(0,1,1),(1,0,1)\}$,
\item $O_{q^2+q+1}:=\{(1,1,1)\}$.
\end{itemize}
%and, using $q$-ary expansion, they can be seen as a partition of the set $\{0,1,q,q+1,q^2,q^2+1,q^2+q,q^2+q+1\}$, which is the set of exponents for the monomials of $F\in \M(x^{q}-x|\F_{q^{3}})$.\footnote{This follows from Theorem \ref{the:charmxq}.}

As discussed in Section \ref{sec:dimxqx}, if we identify a $3$-uple $(a_0,a_1,a_2)$ with the $q$-ary expansion of a number $a_2q^2+a_1q+a_0$, then to each of these orbits we associate the following $\fq$-vector spaces
\begin{itemize}
\item $\V_0:=\fq$,
\item $\V_1:=\{\alpha x+(\alpha x)^q+(\alpha x)^{q^2}: \alpha\in\ff_{q^3}\}$,
\item $\V_{q+1}:=\{\alpha x^{q+1}+\alpha^qx^{q^2+q}+\alpha^{q^2}x^{q^2+1}: \alpha\in\ff_{q^3}\}$,
\item $\V_{q^2+q+1}:=\{\alpha x^{q^2+q+1}: \alpha\in\ff_{q}\}$.
\end{itemize}

Theorem \ref{the:dimxqx} implies that every element of $\M(x^q-x|\ff_{q^3})$ can be uniquely written as a sum of elements in distinct $\V_i$'s. Notice that the sum of the $\fq$-dimensions of the $\V_i$'s is given by $1+3+3+1=2^3$ and  is equal to the dimension of $\M(x^q-x|\ff_{q^3})$.

\bigskip

\item[(2)] Constructing elements in $\M(A|\F_{q^{6}})$ for $A=x^{q^{2}}+x^{q}+x$.

%Let $\mathcal{B}_3$  and $\mathcal{B}_6$ be bases for the $\fq$-vector spaces $\ff_{q^3}$ and $\ff_{q^6}$, respectively. 

%\begin{itemize}
%\item $\{\alpha_1,\alpha_2,\alpha_3\}$, 
%\item $\{\alpha_1 x+\alpha_1^{q} x^{q}, \alpha_2 x+\alpha_2^{q} x^{q},\alpha_3 x+\alpha_3^{q} x^{q},\alpha_4 x+\alpha_4^{q} x^{q},\alpha_5 x+\alpha_5^{q} x^{q},\alpha_6 x+\alpha_6^{q} x^{q}\}$,
%\item $\{\alpha_1 x^{1+q},\alpha_2 x^{1+q},\alpha_3 x^{1+q} \}.
%$
%\end{itemize}
\bigskip

Observe that $A=x^{q^{2}}+x^{q}+x$ is a $q$-additive polynomial dividing $x^{q^{3}}-x$.  By Theorem \ref{the:dimMT}, we can find a $q$-additive polynomial $M$ satisfying $A(M(x))=x^{q^{3}}-x$, namely  $M=x^{q}-x$. 

This polynomial defines a linear map
$$
 \M(x^{q^3}-x|\ff_{q^6})  \xrightarrow{\phi_M} \M(x^{q^{2}}+x^{q}+x|\ff_{q^6})
$$ 
given by $\phi_M(F)=M(F)$.  In other words,  
$$
F^q-F\in\M(x^{q^{2}}+x^{q}+x|\ff_{q^6}),
$$
for all $F\in  \M(x^{q^3}-x|\ff_{q^6})$. Following the strategy of the previous example, we can show that %$\M(x^{q^{3}}-x|\F_{q^{6}})$ is generated by 
$$
\M(x^{q^{3}}-x|\F_{q^{6}})=\ff_{q^3}\oplus\{ \alpha x +(\alpha x)^{q^3}|\alpha \in \ff_{q^6}\}\oplus \{\alpha x^{1+q^3}|\alpha \in\ff_{q^3}\}.
$$
In particular,  $\dim_{\fq}\M(x^{q^{3}}-x|\F_{q^{6}}) =12$, and $\left|\M(x^{q^{3}}-x|\F_{q^{6}})\right|=q^{12}$.
Since $\ker \phi_A=\fq$, we have that 
$$
\left|\M(x^{q^{2}}+x^{q}+x|\ff_{q^6})\right|\geq q^{11}.
$$

\bigskip
 
\item[(3)] An $\mvsp$ different from the ones considered in \cite{carlitz_pol_minimal_61,mills_pol_minimal_64}.

\bigskip

Take $G=x^{q^4+q}-x^{q^3+1}$.  
Notice that $G=M(x^{q^3+1})$ is an element of $\M(x^{q^{2}}+x^{q}+x|\ff_{q^6})$. In particular,  $G$ is an $\mvsp$ with $V_G\subset\ff_{q^3}$ and $\deg G>(q^6-1)/(q^3-1)$. 

Also it is not hard to see that for any $q$-additive polynomial $L$, $G$ is neither a power of $L$ nor of the form $\alpha^{-\deg L} L(\alpha x)+\beta$, for some $\alpha,\beta\in\ff_{q^6}^*$. 

\bigskip
 
\item[(4)] Now assuming that $q$ is odd, we use $\M(x^{q^{2}}+x^{q}+x|\F_{q^{6}})$ to construct elements of 
  $\M(x^{\frac{q^{2}+1}{2}}+x^{\frac{q+1}{2}}+x|\F_{q^{6}})$. 
 
 \bigskip

For $T=x^{\frac{q^{2}+1}{2}}+x^{\frac{q+1}{2}}+x$, we have $x^{q^{2}}+x^{q}+x=T(x^2)/x$. Therefore, Proposition \ref{pro:generalmt} shows that $F^2\in\M(x^{\frac{q^{2}+1}{2}}+x^{\frac{q+1}{2}}+x|\F_{q^{6}})$ for $F\in\M(x^{q^{2}}+x^{q}+x|\F_{q^{6}})$. Using the results obtained in example (2), we conclude that 
$$
\left|\M(x^{\frac{q^{2}+1}{2}}+x^{\frac{q+1}{2}}+x|\F_{q^{6}})\right|\geq (q^{11}-1)/2.
$$
 \end{enumerate}
 
 It is worth mentioning that a computer check for  $q\leq 7$ and $q=2^{n}$ for $n\leq 11$ shows that the lower bound obtained in example $(2)$  is in fact the actual number of $\mvsps$ in the corresponding set $\M(A|\F_{q^{6}})$.
 %\textcolor{red}{Seria isso verdade?}
 
\bibliography{borges_min}{} 
\bibliographystyle{elsarticle-harv}

\end{document}